\documentclass[11pt]{article}

\usepackage{graphicx, float, array, xspace, amsmath, amssymb, amsfonts, bbm}
\usepackage{subcaption}
\usepackage[T1]{fontenc}
\usepackage{setspace}

\usepackage[letterpaper, margin=1in]{geometry}
\numberwithin{equation}{section}

\usepackage[linktocpage=true]{hyperref}
\hypersetup{colorlinks=true, citecolor=blue, linkcolor=blue, urlcolor=blue, pdfauthor={}, pdftitle={},breaklinks=true}

\begin{document}
\pagestyle{empty}      
\begin{center}
\null\vskip0.2in
{\LARGE On the Learnability of Knot Invariants:\\ Representation, Predictability, and Neural Similarity\\[0.5in]}
{Audrey Lindsay$^{a}$, Fabian Ruehle$^{a,b,c}$\\[0.2in]}
{\it 
$^a$ Department of Physics, Northeastern University,\\
Boston, MA 02115, USA\\[3ex]
$^b$ Department of Mathematics, Northeastern University,\\
Boston, MA 02115, USA\\[3ex]
$^c$ NSF Institute for Artificial Intelligence and Fundamental Interactions,\\
Boston, USA
}
\end{center}
\vspace{1cm}
\begin{abstract}
We analyze different aspects of neural network predictions of knot invariants. First, we investigate the impact of different knot representations on the prediction of invariants and find that braid representations work in general the best. Second, we study which knot invariants are easy to learn, with invariants derived from hyperbolic geometry and knot diagrams being very easy to learn, while invariants derived from topological or homological data are harder. Predicting the Arf invariant could not be learned for any representation. Third, we propose a cosine similarity score based on gradient saliency vectors, and a joint misclassification score to uncover similarities in neural networks trained to predict related topological invariants.
\end{abstract}

\clearpage
\pagestyle{empty}
\setstretch{1.0}
\tableofcontents

\setcounter{page}{1}
\pagestyle{plain}
\setstretch{1.0}

\section{Introduction} 
\label{sec:Introduction}
Knot theory is a classical subject in low-dimensional topology. A knot is an embedded circle in a three-dimensional space, typically taken to be $\mathbb{R}^3$ or $S^3$. Under this embedding, two knots are equivalent if they can be related by ambient isotopy. In addition to many rich research areas in pure mathematics that involve knot theory, knots appear in various areas of science, including molecular biology, quantum computing, statistical mechanics, and theoretical physics.

There are many ways to represent a knot. To describe an embedded $S^1$ in $S^3$, one often studies knot diagrams that are obtained by projecting the knot into two dimensions to obtain the knot shadow. To uniquely reconstruct the knot, one keeps track of the over- and under-crossings after projection. Working with the knot projection, Reidemeister showed that two knots are equivalent iff they can be related by a sequence of Reidemeister moves~\cite{Reidemeister:1927aaa}. A key objective in knot theory is to determine whether or not two given knots are equivalent. To this end, multiple knot invariants, i.e., quantities that are invariant under Reidemeister moves (or ambient isotopy), have been proposed. These knot invariants can be geometric or topological in nature, and come in the form of scalar invariants, matrices, polynomials, homological data, etc. 

In practice, some knot invariants are easier to compute than others, and different knot invariants can distinguish different inequivalent knots. Moreover, some knot invariants might be easier or faster to compute in some representation of the knot than in others. There are still many open questions or conjectures concerning knot invariants and their relations. For example, consider the unknot, which is the unique knot with zero crossings. One can now ask which knot invariants uniquely identify the unknot. We know that one of the polynomial invariants called the Alexander polynomial~\cite{Alexander:1928aaa} $A(K)$ cannot distinguish the unknot, since $A(K)=1$ for any knot that is topologically slice. For a different polynomial invariant called the Jones polynomial~\cite{Jones:1985aaa} $J(K)$, it is an open problem whether $J(K)=1$ iff $K$ is the unknot. In general, the Jones polynomial cannot distinguish inequivalent knots, but it might be able to distinguish the unknot. There are categorifications of both polynomials called Heegard Floer and  Khovanov homology, respectively~\cite{Ozsvath:2004aaa, Khovanov:2000aaa}, and both are known to detect the unknot. While we restrict our study to scalar invariants, many of them are related to the above-mentioned polynomials or homology theories, and an interesting question is the extent to which classification algorithms can uncover and make use of these relations.

In terms of computational complexity, computing the Jones polynomial is \#P-hard~\cite{Jaeger:1990aaa}, i.e., unless P=NP, there is no quick way of computing the Jones polynomial. The computational complexity of Khovanov homology is less well-understood, but due to its relation to the Jones polynomial, it is also at least \#P hard. It is another open problem whether the unknot detection problem (independent of the Jones unknotting conjecture) is in P or in NP. The computational complexity and its relation to the representation of the knot was studied in~\cite{BarNatan:2024aaa}. The authors show that determining some invariants, such as the linking number, can be done computationally more efficiently in a 3D representation of the knot than via its 2D diagram. They conjecture that this computational advantage also applies to other knot invariants, such as the hyperbolic knot volume. This motivates the question of which knot representation is best suited for algorithmic implementations of knot invariants.

Given the fundamental nature of knot theory, knots have been fully classified up to a certain number of crossings, and a lot of computer codes and tables have been developed to compute knot invariants. For example, the KnotInfo database~\cite{knotinfo}, which we will be using mostly in this work, contains many topological invariants for knots up 13 crossings, of which there are almost 13000. The software package snappy~\cite{SnapPy} implements and bundles multiple algorithms to compute knot and link invariants. The availability of huge amounts of data and/or code to generate this data, paired with the computational hardness of fundamental problems in knot theory and the abundance of open conjectures and questions about relations of knot invariants, makes it an interesting target for modern machine learning techniques. The first attempt to predict knot invariants using neural networks has been undertaken in~\cite{hughes2016}, while questions such as the unknot detection problem or whether a given knot is ribbon have been addressed with Reinforcement Learning~\cite{Gukov:2020qaj, Gukov:2023kvx}, see~\cite{Gukov:2024buj} for a recent summary. Neural networks were also used by Deepmind to uncover new relations among knot invariants~\cite{Davies:2021Advancing}.

In this paper we compare neural network predictions across different knot representations and knot invariants, as well as relations among neural networks that predict related knot invariants. The results can provide insight into which features of the knot structure are most informative for computing certain knot invariants, which invariants can be computed or approximated efficiently, and whether neural networks can capture relations among knot invariants.

\section{Methods}
\label{sec:Methods}

\subsection{Knot Representations}
\label{sec:Representations}
Given their definition, a perhaps obvious way of representing a knot on a computer is in terms of a closed polygonal curves, or more precisely, in terms of an ordered set of 3D coordinates that describe the embedded knot in 3D. In addition to being well-suited for visualization and physical modeling, this representation can be advantageous for the computation of some knot invariants. However, the computational cost of representing complex knots with 3D coordinates can be high, requiring significant storage and memory. Additionally, while 3D coordinates provide a detailed geometric description, they lack algebraic information and do not capture topological properties of a knot in a straight-forward fashion.

Most other knot representations are based on 2D knot diagrams and encode the crossing information of a knot diagram through sequences of integers, enabling the reconstruction of a knot's planar diagram. In this work, we will study Dowker–Thistlethwaite (DT) codes, Planar Diagram (PD) codes, and Gauss codes. These knot representations are more data-efficient. DT and Gauss codes represent a knot with $N$ crossings using a sequence of $N$ and $2N$ integers, respectively. However, by construction, they encode the crossing information in a non-local fashion, which makes it (at least for humans) harder to picture a knot given in this representation. The PD code provides a more geometrically intuitive encoding, but it is less data-efficient. Indeed, it encodes an $N$-crossing knot in terms of $N$ ordered four-tuples of integers.

Braid words offer an algebraic and efficient way to represent knots and links through the braid group, making them ideal for knot classification and analysis using group theory. They provide a compact representation and are well-suited for computational algorithms. However, they do not offer direct geometric insight into the structure of the knot. Moreover, representing a knot as a (closed) braid requires ambient isotopies and choices of projections that do not allow to represent the knot in terms of its minimal number of crossings. So, while some knots (such as the trefoil with 3 crossings) can be represented by three integers, the braid representation of other knots requires a braid representation that is much larger than $N$, making braids less data-efficient than the other planar knot representations discussed above. For example, the minimal braid word for the 13-crossing knot 13a3143 has length 73. 

While all representations discussed above can be converted into one another, none of them is unique. The 3D representation is only defined up to ambient isotopies, and the representation in terms of polygonal curves involve a random choice of isotopy and a sampling of points along the embedded $S^1$. In practice, the polygonal curve could be represented in different coordinates, using a different number of points, or cyclically permuting all points. In combination to ambient isotopy, the planar invariants involve a choice of projection, which together translate into Reidemeister moves on the planar diagram or Markov moves plus braid relations on the braid word~\cite{Markov:1935aaa}. While there is a canonical choice (i.e., one that leads to the minimal number of crossings), some of the representations are still non-unique. For example, while the four-tuples in the PD code are ordered, there is no natural ordering on the set of four-tuples in the PD code. Similarly, since many generators in the braid group commute, and since some of the Markov moves and braid relations do not change the length of the braid word, there are multiple shortest braid words that represent the same knot, without a canonical choice among them.

\subsection{Knot Invariants}
\label{sec:Dataset}
We focus on scalar knot invariants that are related to the knot diagram, its geometry (for hyperbolic knots), and topological or homological data. In more detail, we look at the following invariants:

\begin{itemize}
    \item \textbf{Alternating (Alt)}: A $\mathbbm{Z}_2$-valued invariant indicating whether the knot admits a planar diagram with alternating over- and under-crossings.
    \item \textbf{Crossing number (Xing)}: An integer-valued invariant given in terms of the minimum number of crossings in any projection of the knot.
    \item \textbf{Unknotting number (Unk)}: An integer-valued invariant that gives minimum number of crossing changes required to transform a knot into the unknot.
    \item \textbf{Signature ($\sigma$)}: An integer-valued invariant, which is computed as the signature of the matrix $M:=V+V^T$, where $V$ is the Seifert matrix. Related to the $s$ and $\tau$ invariants.
    \item \textbf{Determinant (Det)}: An integer-valued invariant corresponding to (the absolute value of) the determinant of $M$ defined above. It is related to the Alexander polynomial; more concretely, it is the absolute value of the Alexander polynomial evaluated at -1.
    \item \textbf{Genus-3D (G3D)}: An integer-valued invariant given by the minimal genus of a Seifert surface embedded in three-dimensional space.
    \item \textbf{Genus-4D (G4D)}: The slice genus is an integer-valued invariant representing the minimal genus of a smoothly embedded surface $\Sigma$ in a 4-ball. $\Sigma$ bounds the knot in $S^3=\partial B_4$. Related to $\sigma$ and the $s$, $\tau$, and Arf invariant.
    \item \textbf{Genus-4D Top. (G4DTop)}: An integer-valued invariant similar to the slice genus but defined for topologically embedded surfaces in the 4-ball.  Related to $\sigma$ and the $s$, $\tau$, and Arf invariant.
    \item \textbf{Rasmussen $\boldsymbol{s}$-invariant ($\boldsymbol{s}$)}: An integer-valued invariant derived from Khovanov homology. It translates to a lower bound on the slice genus. Also related to $\sigma$, $\tau$, and Arf.
    \item \textbf{Ozsv\'{a}th-Szab\'{o} $\boldsymbol{\tau}$-invariant ($\boldsymbol{\tau}$)}: An integer-valued invariant derived from Heegaard Floer homology. It translates to a lower bound on the slice genus. Also related to $\sigma$, $\tau$, and Arf.
    \item \textbf{Arf Invariant (Arf)}: A $\mathbbm{Z}_2$-valued invariant that can be computed from the Seifert matrix. It  is related to the Alexander polynomial and vanishes for knots with slice genus 0. Also related to $s$, $\tau$, $\sigma$.
    \item \textbf{Longitude Length (Lon)}: A real-valued invariant that gives the geometric length of the longitude in a hyperbolic knot complement. This requires the knot to be hyperbolic.
    \item \textbf{Meridian Length (Mer)}:  A real-valued invariant that gives the geometric length of the meridian in a hyperbolic knot complement. This requires the knot to be hyperbolic.
    \item \textbf{Hyperbolic Volume (Vol)}: A real-valued invariant that gives the hyperbolic volume of the knot complement in $S^3$, a fundamental geometric invariant.  This requires the knot to be hyperbolic.
\end{itemize}

While most of this data is present in KnotInfo, we computed invariants that were missing for some of the knots using snappy. In particular, we computed data for longitude length, meridian length, and volume for the 13 crossing knots. In each experiment, we dropped all knots for which the particular invariant we were studying was not known. The results are summarized in Table~\ref{tab:DatasetInfo}. 

The first row shows that for the knots in our dataset, most invariants are known or can be computed with existing code. A notable exception are the topological slice genus (G4DTop) and the Unknotting number (Unk), which are only known for roughly 20 and 40 percent of the knots in the database, respectively. The last three entries are defined for hyperbolic knots only, which are, however, the majority of all prime knots (at least for this range of crossings). 

\begin{table}[t]
\centering
\begin{tabular}{|l|cccccccc|}
\hline
 & Alt & Xing & Unk & $\sigma$ & Det  & G3D & G4D & G4DTop \\
\hline
Data & 12965 & 12965 & 5357 & 12965 & 12965 & 12964 & 12946 & 2973 \\
Class & 2 (2) & 11 (11) & 6 (6) & 21 (11) & 663 (286) & 7 (6) & 7 (7) & 6 (6) \\
\hline
\end{tabular}

\vspace{3mm}

\begin{tabular}{|l|cccccc|}
\hline
 & $s$ & $\tau$ & Arf & Lon & Mer & Vol  \\
\hline
Data & 12965 & 12965 & 12965 & 12955 & 12955 & 12955 \\
Class & 21 (11) & 11 (11) & 2 (2) & $\mathbbm{R}$ & $\mathbbm{R}$ & $\mathbbm{R}$ \\
\hline
\end{tabular}
\caption{Available data for each knot invariant. The unknotting number and the topological slice genus are unknown for many knots. Numbers in brackets specify the number of non-empty classes.}
\label{tab:DatasetInfo}
\end{table}

The limited availability of data can impact the degree and accuracy to which these invariants can be learned as compared to the other invariants with larger training size. Another limitation that can impact training is sparsity of class membership. In the classification task, we set the number of classes to the maximal range of integer values in the dataset. However, this means that for some invariants, many classes have no or only very few members. We give in all cases the number of non-empty classes in the dataset in parentheses in Table~\ref{tab:DatasetInfo}. We see that in particular for the determinant, only about 40 percent of the classes we assign actually have at least one class member, and some determinant values only occur for very few knots. A similar effect can be seen for the signature and the $s$-invariant. In fact, $s$ is always even, and in our dataset $s=2\tau$ (this is false in general but true up to 13 crossing knots). Similarly, $\sigma=-s$ in about 95 percent of the cases in the dataset. We leave it to the NN to learn that the odd classes are not populated, and we find that for the $\sigma$- and $s$-invariant, performance is not impacted. While this data imbalance means that comparison with other invariants is difficult, the relative performance of different knot invariants for these invariants is still meaningful.

For the knot representations, we use zero-padding to make them the same length across different crossing ranges. Since 0 is not used in any of the knot representations other than for padding, this does not alter the knot representation. However, the lengths of the representations differ. DT, Gauss, and PD have lengths $N_\text{max}=13$, $2N_\text{max}=26$ and $4N_\text{max}=52$, respectively. In order to isotope a knot such that it is represented by a braid, one often needs to introduce more crossings. While the braid word length is lower-bounded by $N$, the longest braid word has 73 generators. This means that for example the trefoil has 70 zeros after padding. Instead of zero-padding, one could also perform ambient isotopies (which translate to braid relations and Markov moves) to obtain a braid word of a desired length. For example, we could perform stabilization moves until the braid word has 73 crossings. While we did not perform a systematic study, our results indicate that the way we pad the braid words does not impact the results. Finally, for the coordinate representation, one samples points along the embedded $S^1$ that represents the knot in some isotopy and coordinate system. There is a lot of freedom in how many points one uses to represent the knot as a polygonal curve, with no canonical choice. It is also not clear whether sparse (i.e., the minimum number of points to represent the curve), redundant (include many intermediate points along edges), or dense (sample many points from a smooth embedding) samples are better for the neural network to learn. We choose to represent knots with 174 points (specified by three coordinates each), for a total of 522 real numbers.

\subsection{Neural Networks}
\label{sec:NNs}
We use standard feed-forward neural network architectures for our tests. For the integer invariants, the last layer was chosen for classification (the variance in most integer invariants is low), while we used regression for the real-valued invariants. The networks have three hidden layers with ReLU activation and were implemented in \texttt{pytorch}. 

In order to compare how well-suited different knot representations are for predicting an invariant we want to keep the network expressivity, i.e., the number of parameters, comparable. Since the \texttt{input\_size} of depends on the knot representation, we adapt the input size of the first hidden layer (\texttt{var\_width}) to keep the number of parameters for each knot invariant roughly constant. The width of the subsequent layers is fixed to 100, such that we have
\begin{align}
    \texttt{var\_width}(\texttt{input\_size} + 101) = \text{constant}  
\end{align}
We fix $\texttt{var\_width}=100$ for the coordinate representations that have $\texttt{input\_size}=522$, which means that the constant we use is $62{,}300$. The \texttt{var\_width} of the other knot representation is then scaled up accordingly. Note that the output size (\texttt{num\_classes}) varies for the different invariants but is constant for different representations used to compute the same invariance. 

For training, we use the Adam optimizer with CrossEntropyLoss and MSELoss for classification and regression, respectively, and employ a 90:10 train:test split. We normalize the inputs (and for the regression models, also the outputs) to zero mean and unit variance. Each model is trained for a maximum of 400 epochs, with an early stopping condition if the test loss does not improve for 100 epochs. To determine the best model, we perform a box search over all possible combinations of the following hyperparameters:
\begin{itemize}
    \item \textbf{Learning Rate}: \( \{10^{-5}, 10^{-4}, 10^{-3}, 10^{-2}\} \)
    \item \textbf{Weight Decay}: \( \{10^{-6}, 10^{-4}, 10^{-2}\} \)
    \item \textbf{Batch Size}: \( \{32, 128, 512\}\)
\end{itemize}
The performance of the classifier models is measured by their accuracy, while we use the mean absolute percentage error (MAPE) for the regression models. In more detail, we report $1-\text{MAPE}$ for the unscaled hyperbolic invariants.

\section{Results}
\label{sec:Results}
\subsection{Accuracies}

\begin{figure}[t]
    \centering
    \includegraphics[width=\textwidth]{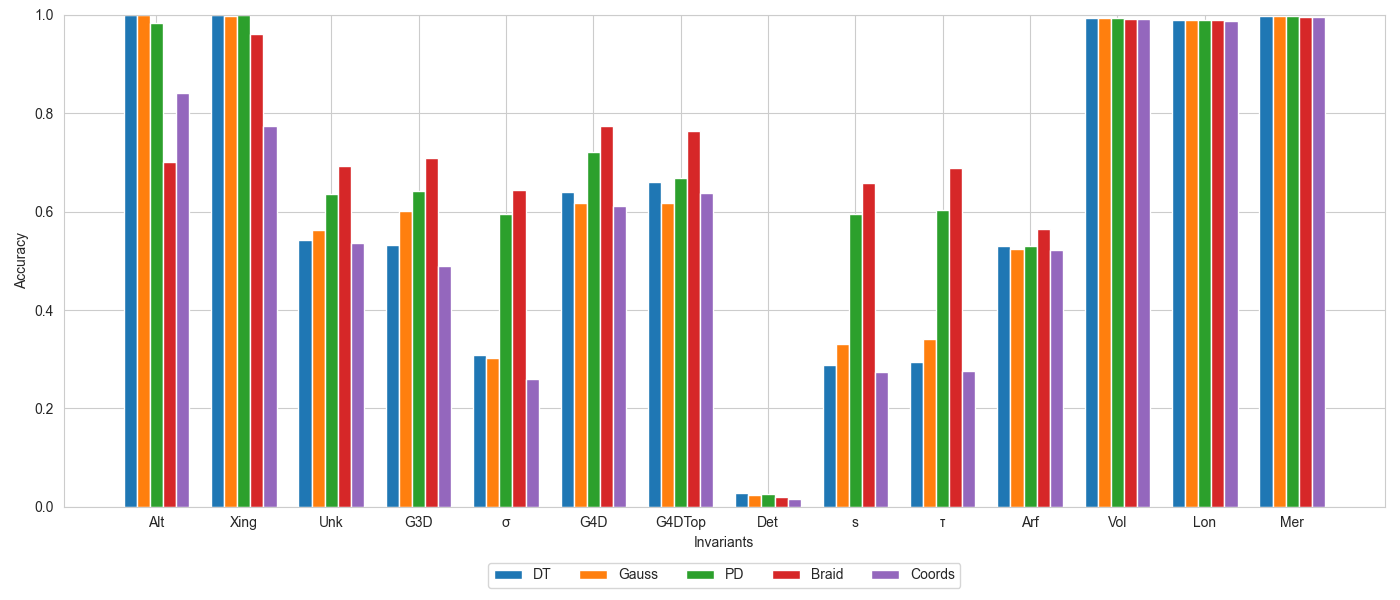}
    \caption{Accuracies of knot representation predictions for different knot invariants.}
    \label{fig:int_invs}
\end{figure}

We present the accuracies of the knot representations in predicting knot invariants in Figure~\ref{fig:int_invs}. For the integer invariants (all but the last 3), we use classification NNs and calculate the accuracy of correct classification. For the last three (hyperbolic) invariants, we calculate the accuracy as 1-MAPE. As a general trend, we find that DT and Gauss codes perform similar, with Gauss codes marginally better. PD codes and braid words outperform DT and Gauss by a wider margin, with braid better than PD. In none of our experiments do the 3D coordinates outperform other representations. Thus, the 3D savings property~\cite{BarNatan:2024aaa}, if applicable to a knot invariant, does not translate into higher prediction accuracy. Of course, the relation between computational savings for a certain representation given some algorithm and the performance of a neural network is not straight-forward.

For the invariants pertaining more directly to the knot diagram, such as whether the knot is alternating and its crossing number, we find that DT, Gauss, and PD perform essentially perfectly. For the crossing number, this is to be expected, since it is encoded straight-forwardly in the features. Indeed, it is simply given by a (fractional) multiple of the nonzero feature entries, and a NN can learn to count these. For braids, the crossing number is upper-bounded by the braid word length, but the actual relation to the crossing number is not straight-forward, making the accuracy of braids in predicting this invariant non-trivial. Similarly, whether a knot is alternating is encoded directly in the sign structure of the DT and Gauss code. It is not too difficult to reconstruct whether a knot is alternating from the PD code, but it requires keeping track of where the repeated indices appear. Again, whether a knot is alternating in braid notation is less trivial to determine, since it requires following the strand around the braid multiple times; in particular, the network has to learn that the ends of the braid are identified. Alternating is a $\mathbbm{Z}_2$-valued quantity and alternating and non-alternating knots are almost evenly distributed in the dataset, meaning a majority classifier would get 50 percent of the cases correct. Finally, the unknotting number is difficult to determine for any knot. Here, PD and Braid notation outperform DT and Gauss, which suggests that the way in which the PD and braid notation encode the crossing information is beneficial for the computation. 

For the topological invariants such as the genera, $s$, $\tau$ and $\sigma$, we see that PD codes and braids have a substantial advantage over DT and Gauss codes. This is an interesting result since it suggests that the homological computation underlying some of these invariants are more readily encoded in braid words than in the other representations. This is perhaps surprising, since Khovanov homology (a variant of which computes $s$) is computed from performing smoothings on crossings, which are more succinctly encoded in the DT, Gauss, or PD codes. As mentioned in Section~\ref{sec:Dataset}, albeit defined and calculated very differently, $\sigma$, $s$ and $\tau$ are linearly related for (almost) all knots in the dataset, 
\begin{align}
    2\tau = s\approx -\sigma\,,
\end{align}
so it is not too surprising that the accuracies for these are the same. These invariants are also related to the slice genus. More precisely, their absolute value gives a lower bound,
\begin{align}
\label{eq:staug4}
g_4 \geq |\tau|\,, \qquad 2 g_4 \geq |s|\,,\qquad 2 g_4 \geq |\sigma|\,. 
\end{align}
The bound is saturated for about 75 percent of the knots in the dataset. We observe that the precision for $\sigma$, $s$ and $\tau$ is about 20 percent (or 15 percent points) lower than that of $g_4$. Also note that $\sigma$, $s$ and $\tau$ are signed invariants (with 11 different values), while $g_4$ is unsigned (with 7 different values). The prediction accuracy for $g_4$, $\sigma$, $s$ and $\tau$ is thus consistent with the NN ``guessing'' $s$ and $\tau$ from $g_4$. We analyze this further in Section~\ref{sec:Similarity}.

Moreover, since a smooth surface bounded by the knot in 4D  can also be used as a surface for the topological genus, and since a surface of a certain genus in 3D could also be used as a surface in 4D, we have
\begin{align}
\label{eq:GeneraBounds}
g_{4,\text{top}} \leq g_4 \leq g_3\,. 
\end{align}
The bound on the left is saturated in almost all (97 percent) of the cases for which both quantities are known in the dataset (which is roughly 23 percent of all knots up to 13 crossings). While this impacts the trustworthiness of a statistical analysis, we find that the accuracies for both four-genera is very similar, except for PD codes scoring better on $g_4$ as compared to $g_{4,\text{top}}$. It would be interesting to study whether this is due to the fact that the PD code needs more training samples or whether this is more systematic. Unfortunately, to the best of our knowledge, there is no scalable method at the moment to compute these invariants, which makes a further study impossible.  The bound on the right-hand side of~\eqref{eq:GeneraBounds} is only saturated in less than 6 percent of the cases, so the NN would have no advantage from using this relation. We observe that the accuracy for $g_3$ is lower than $g_4$, and it would be interesting to understand why predicting $g_3$ is harder than prediction $g_4$.

We also see that none of the representations managed to learn the Arf invariant or the determinant very well. The Arf-invariant is a $\mathbbm{Z}_2$ invariant with roughly equal number of knots having Arf invariant 0 or 1, and thus majority class prediction would lead to an accuracy of 50 percent. For the determinant, as explained above, there is too little data for each individual value of the determinant to expect that NNs can learn this invariant to a high accuracy, which we see reflected in the results. In this case, majority class prediction would give an accuracy of 1 percent. Despite the NNs performing similar to this level, we do not see mode collapse in their output, i.e., their prediction is not just a majority class prediction independent of the input.

For the hyperbolic invariants (volume, meridian length, and longitude length), all representations achieve very high prediction accuracy (within 1 percent of the true value). Since we normalize the inputs and outputs to zero mean and unit variance during training, we report accuracy after transforming the predictions back to the original scale. Without this inverse transformation, the accuracy decreases by approximately 5 percent across all knot representations.

\subsection{Training dynamics}
\begin{figure}[t]
    \centering
    \includegraphics[width=\textwidth]{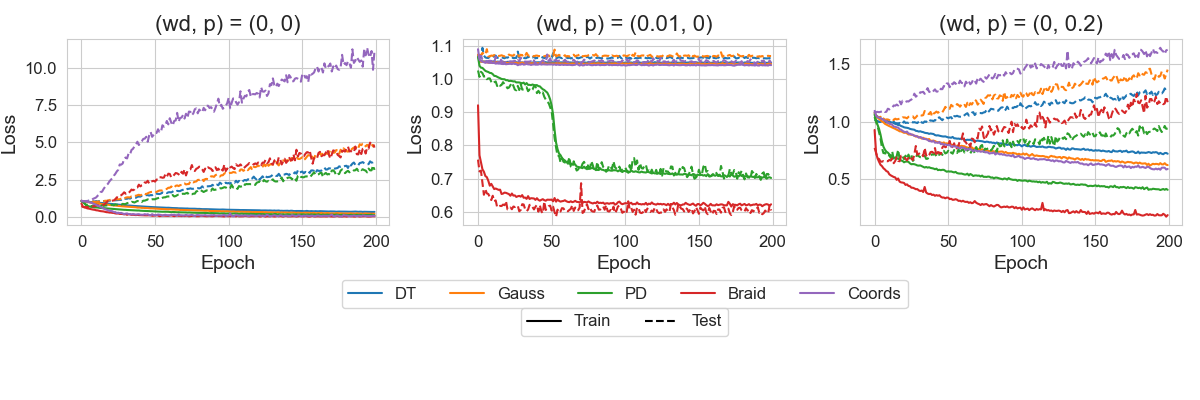}

\caption{Effect of weight decay (wd) and dropout (p) on train and test losses for G4D across different knot representations. Weight decay or early stopping are most effective in preventing over-fitting.}
\label{fig:G4Dloss}
\end{figure}

Across all representations, we found that the models are prone to over-fitting. There are several techniques to ameliorate this:
\begin{itemize}
    \item \textbf{Early stopping:} Stop training before the train loss has converged, or when train and test loss start diverging.
    \item \textbf{Weight decay:} Add weight decay to the optimizer by adding a term $\mathcal{L}_\text{wd}=\lambda\sum|\theta|^2$ to the loss. This pulls the weights towards zero and hence a change in weights needs to be compensated by a corresponding decrease in the loss function.
    \item \textbf{Dropout:} Include dropout layers in the neural network. By randomly dropping nodes, the network cannot fine-tune to a particular signal. 
\end{itemize}
We find all methods to be effective. In particular, we observe that training converges quickly to the optimal model, typically within the first 30 epochs, and for some knot representations and invariants in as little as 10 epochs. This makes early stopping a good candidate, since it means that training time can be very short. Between weight decay and dropout, we find that weight decay is more effective, see Figure~\ref{fig:G4Dloss}. We hence only included early stopping and weight decay in our box search for the hyperparameters. Looking at the test accuracy of the trained model, we find that all methods give the same result.

\subsection{Similarity in NNs across knot invariants}
\label{sec:Similarity}
We want to investigate whether neural networks that are trained to predict different but related knot invariants learn related maps. For this, we study their behavior on the test set. In order to get better statistics, we perform a 70:30 train:test split. This impacts test accuracy, but only by a few percentage points. Since we want to use the same test set for all knot invariants, we remove $g_{4,\text{top}}$ and the unknotting number, since there is only very little data for these. As a measure for similarity across different neural networks, we perform two tests that allow for pair-wise comparison of invariants. We restrict the analysis to braid word representations, since these perform best for almost all invariants.

\begin{figure}[t]
    \begin{subfigure}{0.48\textwidth}
    \centering
    \includegraphics[width=\textwidth]{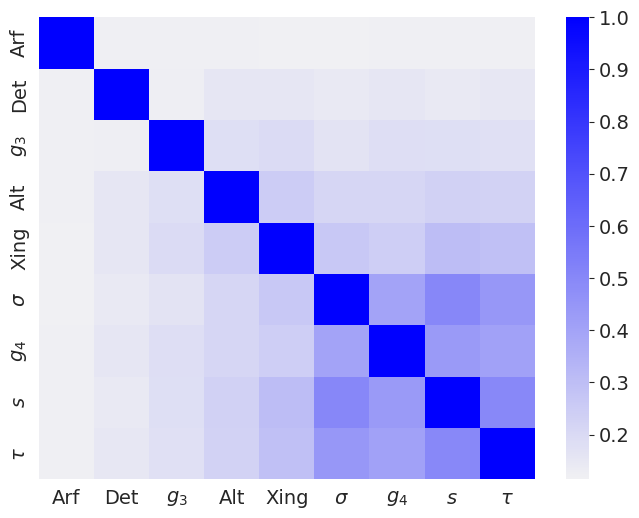}
    \subcaption{Cosine similarity matrix across of gradient saliency scores across knot invariants.}
    \label{fig:CosineSimilarity}
    \end{subfigure}
    \quad
    \begin{subfigure}{0.48\textwidth}
    \centering
    \includegraphics[width=\textwidth]{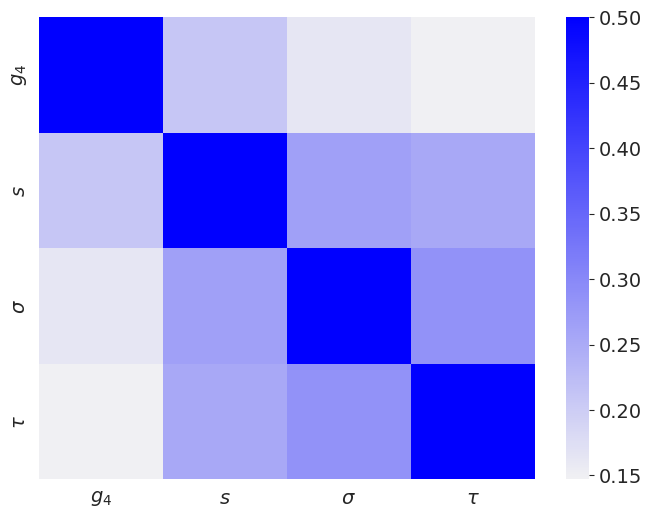}
    \subcaption{Joint misclassification score matrix across knot invariants.}
    \label{fig:MisclassificationScore}
    \end{subfigure}
\caption{Correlation among neural networks for different knot invariants.}
\label{fig:Correlations}
\end{figure}

First, we calculate the gradient saliency $\text{gs}_i(K)$ for each knot $K$ in the test set and each knot invariant $i$,
\begin{align}
    \text{gs}_i(K) = \nabla_K \text{NN}_i(K)\,,
\end{align}
where the gradient is with respect to each letter in the braid word and $\text{NN}_i(K)$ is the NN trained for invariant $i$ and evaluated on $K$. We only compute the saliency for the predicted class, i.e., the output node with the highest probability. Given the sets of gradient vectors $\text{gs}_i$, we compute the pairwise cosine similarity for each knot $K$ between knot invariant $i$ and $j$,
\begin{align}
    \text{c}_{i,j}(K) = \frac{\text{gs}_{i}(K) \cdot \text{gs}_{j}(K)}{||\text{gs}_{i}(K)|| \text{gs}_{j}(K)||}\,.
\end{align}
From this, we compute a scalar summary statistic by averaging over the cosine similarity of all knots,
\begin{align}
    \overline{\text{c}}_{i,j} = \frac{1}{N}\sum_{K}\text{c}_{i,j}(K)\,,
\end{align}
where $N$ is the number of knots $K$ in the test set. A heatmap of the matrix $\overline{\text{c}}_{i,j}$ is given in Figure~\ref{fig:CosineSimilarity}.

Second, we determine the set $S_i$ of knots that are misclassified for each knot invariant $i$. From these sets, we compute a joint misclassification score
\begin{align}
    \text{miss}_{i,j} = \frac{|S_i\cap S_j|}{|S_i|+|S_j|}\,,
\end{align}
where $|\cdot|$ denotes the cardinality of the set. With this definition $0\leq \text{miss}_{i,j} \leq .5$. We focus only on those knot invariants that showed substantial overlap in the cosine similarity measure. We plot a heatmap of the misclassification score in Figure~\ref{fig:MisclassificationScore}.

From either measure, we observe a strong correlation between $g_4$, $\sigma$, $s$, and $\tau$. The Arf invariant and Det are not strong correlated with any other knot invariant, which makes sense since the NN has essentially failed to learn these invariants and is hence not making reliable predictions. The crossing number and whether the knot is alternating also show some correlation. Moreover, $\sigma$, $s$ and $\tau$ show a very strong correlation, and all are also strongly correlated with $g_4$. Hence these measures correctly identifies related but distinct invariants.

\subsection{Dependence on hyperparameters and architecture}
The networks in our study have around 70k parameters for a dataset size of 13k. We also ran the same experiments with larger NNs (with around 2M parameters), as well as with different activation functions, but the accuracies did not improve. This indicates that the accuracy is limited by training size and not by model size. While the hyperparameter box search showed significant variation in the final accuracies, the highest accuracies were obtained for hyperparameters at or close to the pytorch defaults.

The poor performance of the coordinate representation begs the question of whether the performance can be improved. One possibility would be an ablation study in the number of sampled points on the curve. Another possibility would be using equivariant neural networks that take into account some of the symmetries of $\mathbbm{R}^3$~\cite{Geiger:2022aaa}. Of course, a fully equivariant neural network with respect to ambient isotopy is out of reach. As a perhaps simpler approach, one can study how neural network architectures targeted towards vision tasks perform for the coordinate representation. While a full-fledged study of this is beyond the scope of this paper, we coded a simple ResNet~\cite{He:2015ResNet} with three residual blocks, consisting of batch norms, convolutional layers and skip connections, followed by pooling and fully connected layers with a total of 170k parameters. We choose kernel sizes of height 3 to keep the $(x,y,z)$-coordinates of the points grouped together. We did not perform a systematic hyperparameter search for these architectures but found a modest improvement of the accuracies for the crossing number and for alternating to 80 and 96 percent, respectively. The other invariants did not improve statistically significantly. It is interesting that the invariants related more directly to the knot diagram improved with this architecture change. 

Another question is whether the poor performance on predicting the Arf invariant can be improved. Motivated by the similarity of braid words to natural language, we tried a transformer~\cite{Vaswani:2017aaa} for predicting the Arf invariant. The transformer was shown to give slightly better performance on braid words than FFNN for the unknot decision problem in~\cite{Gukov:2020qaj}. We used a simple transformer as encoder with 8 attention heads and 12 layers (around 2.5M parameters) and trained it for 500k steps with batch-size 64, Adam optimizer with learning rate $10^{-4}$, 8 attention heads, 12 layers, and 128-dimensional embedding layer with sinusoidal positional encoding.  We also tried the sequence to sequence model~\cite{Int2Int} with a decoder with similar architecture to the encoder. Without proper hyperparameter tuning, the accuracy for the Arf invariant improved slightly to 62 percent, which is better than the FFNN, but still lacking behind the other invariants. A problem could be that the training set with is too small for a transformer, and it would be interesting to see whether the performance can be improved with a larger dataset. We leave this for future study.

\section{Conclusion}
\label{sec:Conclusions}
In our study, we compare how well different knot representations allow for computing a given knot invariant, as well as how it is to predict a particular knot invariant as compared to other invariants. For both, we observe large variances.

For the choice of representations, braid words perform in general the best, while 3D coordinates, DT, and PD codes perform the worst. Comparing the prediction accuracy across invariants, we find that diagram invariants and hyperbolic invariants are the easiest to predict (although their accuracy is computed differently), while the topological invariants are harder. In particular, we find that none of the architectures or knot representations allowed to learn the Arf invariant. In contrast, the 3D and 4D genus can be predicted from the braid word with almost 80 percent accuracy, and related quantities such as the signature and the $\tau$- and $s$-invariant can be learned with about 70 percent accuracy. The determinant and the topological 4-genus suffer from a lack of data to make reliable predictions.

The best hyperparameters selected by the box search do not significantly improve the test accuracy over the standard \texttt{pytorch} parameters. Via an ablation study in model size, we find that accuracy is limited by training size and not model complexity. Using a different architecture did not significantly affect the accuracies. During training, all models converged quickly but showed signs of overfitting, which was effectively prevented by early stopping and weight decay.

Some of the invariants in the dataset are related, although they are computed in very different ways. In particular, the values of $s$, $\tau$ and $\sigma$ are strongly related with one another, and, to a lesser extent, also with $g_4$. To quantify how similarly neural networks trained to predict these invariants learn their mappings, we propose two measures. First, we compute the gradient saliency of the networks and then compute a cosine similarity score for these gradients. Second, we look at joint misclassifications across different neural networks and invariants. These methods successfully discover different but related topological invariants in the dataset. 

In the future, it would be interesting to extend the similarity study of neural networks to non-scalar invariants such as the Jones or Alexander polynomial, or Khovanov and Heegard Floer homology. Comparing scores for $s$ and $\tau$ with the corresponding homological data could provide further insight. In particular, it would be interesting to understand why the Arf invariant is so elusive, and whether the good performance of braid representations can be leveraged to implement more efficient algorithms for computing knot invariants. Moreover, examining cases where different invariants lead to divergent predictions may reveal how the neural network learns distinct mappings in these instances. To understand the (average) computational complexity or approximation schemes for the different knot invariants, it would also be useful to study the maximum accuracy that can be obtained for the invariants. To this end, a systematic study of other architectures, such as ResNets or transformers that we started to look at would be helpful.

\section*{Acknowledgments}
The work of FR is supported by the NSF grants PHY-2210333 and PHY-2019786 (The NSF AI Institute for Artificial Intelligence and Fundamental Interactions), as well as startup funding from Northeastern University. The computations were carried out on Northeastern's Discovery cluster. We are grateful for support from CMSA for the Mathematics and Machine Learning Program during the Fall 2024 semester, where this project was initiated.

\bibliographystyle{JHEP}
\bibliography{bibliography}
\end{document}